\input amstex
\documentstyle{amsppt}
\magnification=1200
\hcorrection{.25in}
\advance\vsize-.75in

\catcode`\@=11
\def\logo\@{}
\catcode`\@=\active  
\NoBlackBoxes
\topmatter
\title   Asymptotics via Steepest Descent for an Operator Riemann-Hilbert Problem\endtitle
\author Spyridon Kamvissis \endauthor
\endtopmatter
\document

\define\[{\left[}%
\define\]{\right]}%
\define\({\left(}%
\define\){\right)}%

\bigskip

ABSTRACT

\bigskip

In this paper, we take the first step towards an extension of  
the nonlinear steepest descent method
of Deift, Its  and Zhou  to the case of operator Riemann-Hilbert problems.
In particular, we provide long range asymptotics for a Fredholm determinant 
arising in the computation of the probability of finding a
string of n adjacent parallel spins up in the antiferromagnetic ground state of
the spin 1/2 XXX Heisenberg Chain. Such a determinant can be expressed in 
terms of the solution of an  operator Riemann-Hilbert factorization problem.

\newpage

1. INTRODUCTION-MOTIVATION

\bigskip

It has been established by Korepin, Izergin and Bogoliubov [KIB] that many 
important problems arising in the theory of Quantum Inverse Scattering 
can be reduced to Classical Inverse Scattering Problems.  In particular, 
the correlation functions for quantum solvable models can be described by
solutions of classical 
integrable equations, or, equivalently, are reducible to the solution of 
Riemann-Hilbert problems. An important observation is that many such 
quantities are described by Fredholm determinants of particular (`integrable') 
integral operators. As shown by different authors 
(see [DIZ] for a concise proof) the computation of these determinants 
reduces to the solution of a Riemann-Hilbert problem. Long range or long 
time asymptotics for such a quantity  can then be recovered 
via an asymptotic analysis of an oscillatory Riemann-Hilbert 
problem.\footnote{See e.g.[DZ1] and [DIZ] for the Ising Chain and the XXO 
model. Such analysis has also been performed in the context of 
long time asymptotics for soliton equations by several authors; see
e.g. [DZ2], [K1], [DVZ], [DKKZ], [K2].}  
In general, such problems are not local; in fact they can be considered as 
operator Riemann-Hilbert problems, which makes their analysis harder. 
In this paper we present a solution of one of these problems, by
reducing  it to two distinct (and essentially different)
matrix (2x2) problems and thus easily extracting asymptotics.

We are motivated by the following observation [KIEU].

FACT 1.1. Let $\Pi(n,\psi, \phi)$ be the 
probability of finding a string of 
n adjacent parallel spins up in the
antiferromagnetic ground state of the spin 1/2 XXX Heisenberg chain. 
Here $\psi$ is an angle related to
the external magnetic field $h,$\footnote{In fact   $\psi$ is defined by
$e^{-i\psi}= -i {{e^{-2L}-i}\over{e^{-2L}+i}}, -\pi <\psi<0, $ where
$cosh2L= 2/h, L>0$.} and $\phi(z)$ is 
a dual quantum field, which can be considered as a holomorphic
function in a neighborhood of the unit circle.
Then, $\Pi(n,\psi,\phi)=  {{ (0|det (1+V) |0) } \over  {det (1+K) }}$,
where $0$ is the Fock vacuum, $K$ is the integral operator 
on $L^2 [0, \infty) $ with kernel
$1 \over {\pi (1+(\lambda-\mu)^2) }$ and $V$ is an integral operator defined
below.

Our aim is the following.

THEOREM 1.2. Let $\psi$ be as above and $\phi (z) $ be an entire function.
Let $P(n, \psi, \phi)= det (I+V)$, where
$V$ is given by (2.1) below. As $n \to \infty$, we have
$$
\aligned
P(n, \psi, \phi) \sim const~exp (n^2 log|sin(\psi/2)| + O(logn)).
\endaligned
\tag1.1
$$
Note that there is no $\phi$ dependence.

REMARK. Theorem 1.2 falls short of computing the asymptotics for the
actual correlation function $\Pi (n, \psi, \phi)$. This may require the
complete asymptotic expansion of $P (n, \psi, \phi)$. We plan to consider this
problem in a later publication.

In the next section,  we set up the operator
Riemann-Hilbert problem that enables us to compute the Fredholm determinant
$P(n, \psi, \phi)$. Using deformations inspired by [DZ2], [FIK] and [K2]
we will reduce this problem to  two standard matrix ones, which we can solve
explicitly (as $n \to \infty$).

\bigskip

2. THE OPERATOR RIEMANN-HILBERT PROBLEM

\bigskip

We are interested in the long range asymptotics of the
Fredholm determinant $det (I + V^{(n)})|_{\gamma=1}$ where the integral 
operator $V^{(n)}$ acts on  a function $f$ as follows:
$$
\aligned
(V^{(n)} f) (z_1) = \int_C V^{(n)} (z_1, z_2) f(z_2) dz_2,
\endaligned
\tag2.1
$$
where $C$ is the contour $\theta \to z= exp (i\theta), ~-\psi<\theta<2\pi+\psi$
and
$$
\aligned
V^{(n)} (z_1, z_2) = -\gamma {i \over {2\pi}} { 1 \over {z_1-z_2}}
[e_+^{(n)} (z_1) e_-^{(n)} (z_2) r(z_1, z_2) 
-e_-^{(n)} (z_1)e_+^{(n)}(z_2) r(z_2, z_1) ],\\
r(z_1, z_2) = {{ 2(z_2-1) (z_1-1) } \over
{2 (z_2-1) (z_1-1) + z_1 -z_2 }}, ~~~~
e_{\pm}^{(n)} (z) = z^{\mp n/2} e^{\pm \phi(z)/2}.
\endaligned
\tag2.2
$$
Here $\psi \in (-\pi,0)$ is defined as in Theorem 1.1 and is related to an
external magnetic field.

It is established that for special 
(`integrable') operators where the kernel is of the form
$$
\aligned
V(z_1, z_2) = \Sigma_1^N {{a_j(z_1) b_j(z_2)} \over {z_1-z_2}}, ~~with~~~
\Sigma_1^N a_j(z) b_j(z) =0,
\endaligned
\tag2.3
$$
the computation of the related Fredholm determinant can be reduced to the 
solution of an NxN matrix Riemann-Hilbert problem (see [DIZ] for the proof).
In our case, the kernel (2.2) is not of the desired type. However, using
$$
\aligned
r(z_1, z_2) = \int_0^{\infty} 
e^{-s+ {s \over 2}{{z_1+1}\over{z_1-1}} -{s \over 2} {{z_2+1}\over { z_2-1}} } ds
\endaligned
\tag2.4
$$
we obtain ([FIK])
$$
\aligned
V^{(n)} (z_1, z_2) = -\gamma {i \over {2\pi}} \int_0^{\infty}
e_+^{(n)} (z_1|s) e_-^{(n)} (z_2|s) - e_-^{(n)} (z_1|s) e_+^{(n)} (z_2|s) ds,
\endaligned
\tag2.5
$$
where
$$
\aligned
e_{\pm}^{(n)} (z|s) = (z^{-n} exp (\phi(z) + s {{z+1} \over{z-1}} ) )^{\pm1/2}
e^{-s/2} .
\endaligned
\tag2.6
$$
In other words, we end up with an operator Riemann-Hilbert problem as 
follows.

THEOREM 2.1. [FIK] Let $\Psi (z)$  (similarly $M(z)$)
be a 2x2 matrix of integral operators
depending on the complex parameter $z$: 
$$
\aligned
(\Psi_{jk}(z) f) (s) = \int_0^{\infty} \Psi_{jk} (z|s,t) f(t) dy,~~~
f \in L^2[0,\infty),
\endaligned
\tag2.7
$$
such that

1. $\Psi (z)$ is analytic for $z \in \Bbb C \setminus \bar C$.

2. $\Psi (\infty) = I \pmatrix &1 &0 \\ &0 &1 \endpmatrix,$ where
$I(s,t) = \delta (s-t).$

3. $\Psi_+(z) = \Psi_- (z)  M(z), ~~~~~z \in C$,
where $\Psi_+$ and $\Psi_-$ are the normal limits of $\Psi$ from inside
and outside the unit circle respectively, and
$$
\aligned
M(z|s,t) = \pmatrix &\delta (s-t) + \gamma e_-^{(n)} (z|s) e_+^{(n)} (z|t)
&-\gamma e_+^{(n)} (z|s) e_+^{(n)} (z|t) \\ 
&\gamma e_-^{(n)} (z|s) e_-^{(n)} (z|t) 
&\delta (s-t) - \gamma e_+^{(n)} (z|s) e_-^{(n)} (z|t) \endpmatrix,
\endaligned
\tag2.8a
$$
or equivalently,
$$
\aligned
M(z) = \pmatrix &I+ \gamma P^T & -Q(z) z^{-n} e^{\phi(z)} \\
&Q(1/z) z^n e^{-\phi(z)} & I- \gamma P(z) \endpmatrix,~~~z \in C,
\endaligned
\tag2.8b
$$
where 
$P$ and $Q$ are given by the kernels
$$
\aligned
P(z|s,t) = exp [{1 \over 2} (s-t) {{z+1} \over {z-1}} - {{s+t} \over 2} ],\\
Q(z|s,t) = exp [{1\over 2} (s+t) {{z+1} \over {z-1}} - {{s+t} \over 2} ].
\endaligned
\tag2.8c
$$
Then, $P(n, \psi, \phi) = det(I+V^{(n)})$ satisfies the relation
$$
\aligned
{{P(n+1, \psi, \phi)}\over {P(n, \psi, \phi)}} = det \Psi_{22} (0).
\endaligned
\tag2.9
$$

PROOF: The proof is analogous to the proof of the pure matrix equivalent
(see [DIZ]); it is essentially given in [FIK].

REMARK: We note that $P$ and $I-P$ are orthogonal 
projection operators. In fact, $P^* =P$ and
$P^2 (z) =P(z)$; also $P(z)Q(z) = Q(z),$
and $Q(z) Q(1/z) = P(z)$.
By definition, $P^T (z|s,t) = P(z|t,s)$. 
Formulae (2.8c) define  $Q$ at least  for $Rez<1$. $P$ is
defined for at least  the unit circle minus the point $z=1$.
An essential singularity exists at $z=1$.

From now on, we focus our attention to the (physically interesting case)
$\gamma =1$.

We use the orthogonal decomposition into
$Im(I-P)$ and $ImP$. Our Riemann-Hilbert problem defined by (2.8)
splits into two distinct ones, of different nature. For convenience,
we consider the jump contour as the whole unit circle, with jump
$$
\aligned
M^* = M, ~~~z \in C,\\
M^*= I, ~~~z \in (|z|=1) \setminus \bar C.
\endaligned
$$
We then have 
$$
\aligned
(\Psi (I-P) )_+ = \Psi_- (M^*(I-P)),\\
(\Psi P )_+ = \Psi_- M^*P,
\endaligned
$$
across the unit circle, with identity asymptotics at infinity.
In other words, we have two new problems. Let $\Phi = \Psi (I-P)$
inside the unit circle and $\Phi = \Psi$ outside the unit circle. Then
$$
\aligned
\Phi_+ = \Phi_- M^* (I-P),~~~z \in C,\\
lim_{z \to \infty} \Phi =I.
\endaligned
\tag2.10
$$
Furthermore, $\Phi_+$ has an essential singularity at $z=1$.

Likewise, let $\chi = \Psi P$ inside the unit circle and
$\chi = \Psi$ outside the unit circle. Then
$$
\aligned
\chi_+ = \chi_- M^* P, ~~~z \in C,\\
lim_{z \to \infty} \chi = I.
\endaligned
\tag2.11
$$
Again, $\chi_+$ is essentially singular at $z=1$.

Note that 
$$
\aligned
\Psi (0) = \chi (0) + \Phi (0).
\endaligned
\tag2.12
$$
These new problems are degenerate since none of $P, ~ I-P$ is invertible.
They do become nondegenerate once we restrict $\Phi(z), \chi(z)~~~|z|<1$ and
$M^*(z),~~|z|=1$ to  the two orthogonal subspaces $I-P $ and $P$ respectively. 
Note that for
each of the new problems the jump across $\{|z|=1\} \setminus \bar C$ is again
the identity. Even more, the singularity at $z=1$ is removed.

On $Im(I-P)$ we have $(I-P)|_{Im(I-P)} =I.$ Thus,\footnote{Note here that
$Im~M(I-P)$ is independent of $z$ by continuity.}
$$
\aligned
\Phi|_{Im(I-P),+} = \Phi|_{Im M(I-P),-} M|_{Im(I-P)},~~~~~where\\
M|_{Im(I-P)} (z) = \pmatrix &I & -Q(z)z^{-n} e^{\phi(z)} \\
&Q(1/z) z^n e^{-\phi(z)} &I \endpmatrix,~~~~z \in C.
\endaligned
\tag2.13
$$
This is essentially the problem appearing in the case $0<\gamma<1$,
treated in [FIK]. We have the obvious factorization 
(note $Q(z) Q(1/z)_{|Im(I-P)} =0$)
$$
\aligned
M|_{Im(I-P)} (z) = M_U(z) M_L^{-1} (z),~~~~~where\\
M_U(z) =\pmatrix &I &-Q(z) z^{-n} e^{\phi(z)} \\
&0 &I \endpmatrix,\\
M_L(z)= \pmatrix &I &0 \\
&-Q(1/z) z^n e^{-\phi(z)} &I \endpmatrix.
\endaligned
$$
By using the analyticity of $Q(z)$ for $|z|<1$ and $Q(1/z)$ for 
$|z|>1$, we
deform the problem as follows.

Let $\Sigma$ be the augmented contour consisting of the union of:

1. The contour $C$.

2. A smooth curve $C_{int}$ joining the endpoints of the contour $C$, lying
entirely within the open disc $|z|<1$ and close to $\{|z|=1\} \setminus C$.

3. A smooth curve $C_{ext}$ joining the endpoints of the contour $C$, lying
entirely in $|z|<1$.

All contours are meant to have a counterclockwise orientation.
The complex plane is now divided into three regions.

1. The region containing $0$, say $R_1$.

2. The unbounded region, say $R_2$.

3. The region between $C_{int}$ and $C_{ext}$, say $R_3$.

\bigskip

Let 
$$
\aligned
\tilde \Phi = \Phi|_{Im(I-P)} M_L,~~~~z \in R_1,\\
= \Phi|_{Im M(I-P)} M_U,~~~~z \in R_2,\\
= \Phi|_{Im(I-P)} = \Phi|_{Im M(I-P)},~~~ z \in R_3.
\endaligned
$$
Then, the jumps for $\tilde \Phi$ are equal to $M_U$ on $C_{ext}$ and
$M_L^{-1}$ on $C_{int}$. On $C$, there is no jump since
$\Psi_+ M_{L+} = \Psi_- M_{U-}$. As $n \to \infty$, we trivially
get the identity solution. So,
$$\aligned
\chi_{|Im(I-P)} (0) \sim I,~~~~~as~~~~n \to \infty.
\endaligned
\tag 2.14
$$
Note that the essential singularity at $z=1$ never plays a role since
$z=1$ is not on $C$.

Also, clearly $\chi_{|ImP} (0) =0, \Phi_{|Im(I-P)} (0) = 0.$
So, it remains to calculate $\Phi_{|ImP} (0)$; in fact,
$$
\aligned
det \Psi_{22} (0) \sim det \Phi|_{ImP;22} (0), ~~~as~~~n \to \infty.
\endaligned
\tag2.15
$$
Hence, we can simply focus our attention
on $ImP$.

On $ImP$, $Q(z)={{1-z} \over 2}I.$ In particular, no
essential singularity exists and  $Q(z)$ and $Q(1/z)$
can be defined meromorphically
in the whole complex plane using $Q(z) Q(1/z) =I$.
Our operator Riemann-Hilbert problem has jump matrix
$$
\aligned
M|_{ImP} (z) = I \pmatrix  &2  &- Q(z)~z^{-n} e^{\phi(z)} \\
& Q(1/z)~z^n e^{-\phi(z)} &0 \endpmatrix,~~~z \in C.
\endaligned
\tag2.16
$$
The conjugation used in the previous  case is useless as the resulting 
scalar problem has jump with determinant 0. Recognizing that we have
here an operator version of a `shock'-type Riemann-Hilbert problem, we
have to use a conjugation with a function defined appropriately as a 
radical (cf. [DVZ], [K2], [DIZ]).

Let $ \alpha = -sin^2{\psi/2} <0.$ As in [DIZ], p.218, we define the function
$g$ such that
$$
\aligned
g(z)~~is~~~analytic~~~in~~~\Bbb C \setminus \bar C;\\
g(z) \to 1,~~as~~~z~~~\to \infty;\\
g_+ (z) g_-(z) = {\alpha \over z},~~z~~\in ~~~C;\\
|{g_+ \over g_-}|<1, ~~z~~~\in~~~C.
\endaligned
$$
It actually follows that 
$$
\aligned
g(z) = {{ ((z-e^{-i\psi} ) (z-e^{i\psi}) )^{1/2} + z-1} \over {2z}},~~~so~~
g(0) = sin^2{\psi \over 2}.
\endaligned
\tag2.17
$$

We now define a new operator valued 2x2 matrix by
$$
\aligned
F_n(z) = \alpha^{n\sigma_3/2} \Phi_{|ImP}(z) 
g(z)^{n\sigma_3} \alpha^{-n\sigma_3/2},
~~|z| <1, \\
F_n (z) = \alpha^{n\sigma_3/2} \Phi|_{Im MP}
(z) g(z)^{n\sigma_3} 
\alpha^{-n\sigma_3/2}, ~~~~|z|>1,
\endaligned
\tag2.18
$$
where $\sigma_3 = diag (1, -1)$ is a Pauli matrix. Then, $F_n(\infty)=I$ and
the jump matrix for $F_n $ is 
$$
\aligned
M^{ImP}_{g,n} (z)=  
\pmatrix & 2 ( {g_+ \over g_-})^n  I  &-Q(z)e^{\phi(z)} \\
&Q(1/z) e^{-\phi(z)} & 0 \endpmatrix, ~~~z \in C.
\endaligned 
\tag2.19
$$
As $n \to \infty$,
$$
\aligned
M^{ImP}_{g,n} (z) \to \tilde M = I \pmatrix 
&0 &- {{1-z}\over2}e^{\phi(z)} \\&  {2 \over {1-z}} e ^{-\phi(z)} &0 \endpmatrix, ~~~z \in C.
\endaligned
\tag2.20
$$
It is thus appropriate to consider the Riemann-Hilbert
$$
\aligned
\tilde F_+ = \tilde F_- \tilde M,~~~z \in C,\\
\tilde F(\infty) =I.
\endaligned
\tag2.21
$$
hoping that $\tilde F= F^{\infty} = lim_{n \to \infty} F_n $.
It turns out that the solution of (2.21) has a pole at $z=1$.
This  problem  can be conjugated to one with jump 
equal to $I \pmatrix &i &0 \\ &0 &-i \endpmatrix$
and the solution can be derived
easily. We diagonalize $\tilde M$ as follows.
$$
\aligned
\tilde M = S D S^{-1},\\
S= \pmatrix & 1 &1 \\
&{{2i} \over{z-1}} e^{-\phi(z)} &{{-2i} \over {z-1}} e^{-\phi(z)} \endpmatrix,\\
D= diag (i, -i).
\endaligned
\tag2.22
$$
Now $D$ can be factorized as 
$$
\aligned
D=B_-B_+^{-1},\\ B= diag(\beta^{-1}, \beta),
\endaligned
\tag2.23
$$
where
$$
\aligned
\beta(z) = ( {{z-e^{-i\psi} } \over {z-e^{i\psi}} })^{1/4},\\ \beta(\infty)=1.
\endaligned
\tag2.24
$$
Hence $(\tilde F S B)_+ = (\tilde F S B)_- $ across the jump contour;
but $\tilde F S B  \sim S$  near $z=\infty$, hence in fact
$\tilde F S B = S$ and 
$$
\aligned
\tilde F = S B^{-1} S^{-1} = {I\over 2} \pmatrix
&\beta +\beta^{-1} &{i(z-1) \over 4} e^\phi (\beta- \beta^{-1}) \\
&{i \over {z-1}} (\beta-\beta^{-1}) &\beta + \beta^{-1} \endpmatrix.
\endaligned
\tag2.25
$$
Even though  $\tilde F$ has a pole at $z=1$, the (22)-entry which concerns
us does not. It is thus possible to show that $F_{n,22} \to \tilde F_{22}$,
as $n \to \infty$.
The following formulae
for the Fredholm determinants below follow immediately.
$$
\aligned
det \tilde F_{22} (0) = sin {\psi \over 2},\\
det \Psi_{22}(0) = -sin^{2n+1} {\psi \over 2}.
\endaligned
\tag2.26
$$
Using  (2.9), (2.17) and (2.18) we get
$$
\aligned
{P_{n+1} \over P_n} = sin^{2n+1} {\psi \over 2} ~~(1+0({logn \over n})).
\endaligned
$$
Theorem 1.2 now follows.

PROOF OF THEOREM 1.2. A rigorous proof that $F^{\infty}_{22} = \tilde F_{22}$
(as well as formula (2.14)) 
requires a Beals-Coifman type formula
which reduces the  operator Riemann-Hilbert problem to
a singular integral equation. This can be easily done and the details of
the proof are essentially contained in [DIZ], p.159. The Cauchy
operator involved takes bounded operators in $L^2$ to 
bounded operators in $L^2$ and is itself bounded.

A small complication arises because the
convergence of $|{g_+ \over g_-}|^n$ is not uniform near the endpoints of
$C$. This is dealt with  by constructing a parametrix near those
points, as suggested in [DIZ] for the analogous matrix problem. $g$ being
scalar, there is no complication due to the operator nature of the
underlying problem. 

Finally,  one may worry about the fact that
$\tilde F$ has a pole at $z=1$ and the possible effect on the validity of the
limiting procedure. This, however has been shown not to be the case,
since any meromorphic problem can be made into a holomorphic problem with
an extra jump on a circle around the pole (see e.g.[DKKZ]). The important fact
is that the (22)-entry of $\tilde F$ (which is all
we are interested in) has no singularity.

The limit $F^{\infty}_{22}$ is meant to be taken in the 
trace class norm. The result for determinants follows readily and the proof
of Theorem 1.2 is complete.

\bigskip

REMARKS:

1. It may initially seem  that an essential singularity at $z=1$ is introduced in the
Riemann-Hilbert problem when one multiples the unknown with the operators $P$
and $I-P$ (cf. remark after  the proof of theorem 2.1). However, this 
singularity vanishes once we restrict our operators to smaller subspaces.

\bigskip

2.  In Theorem 1.2, we have been assuming that 
$\phi(z)$ has an analytic extension
on the complex plane. This is not necessarily true in the original physical
problem (see Fact 1.1); however, we can always
approximate $\phi$ by an analytic $\tilde \phi$ that enables our calculations
to go through, and recover the same result (see [DZ2] for details). 

\bigskip

3.  In conclusion, we would like to point out that the importance of
the above procedure lies in the fact that operator Riemann-Hilbert problems
appear in several contexts, apart from the computation of
correlation functions for exactly solvable models of statistical
mechanics. For example, it is well known that the
inverse problem for integrable equations in 2+1 dimensions (Davey-Stewartson,
KP, etc.) can often be expressed as an operator Riemann-Hilbert problem, and
hence deformed along the lines described here.

\bigskip

4.  We also note that, even though operators do not generally commute,
and hence even `scalar' operator Riemann-Hilbert problems are not trivial
to solve explicitly, scalar conjugating functions (like $g$ of (2.10)) can
still play the important `deforming' role they play in the standard
matrix factorization problems.

\bigskip

ACKNOWLEDGEMENTS. We thank Nicos Papanicolaou for pointing our attention to
this problem and for his hospitality at the Physics Department of the
University of Crete. We also thank V.E. Korepin and A.R.Its 
for several important comments.
A previous version of the present paper has appeared  in the
preprints series of the Max-Planck Institute for Mathematics
in the Sciences, Leipzig, n.1998-38.

\bigskip

BIBLIOGRAPHY

\bigskip

[DIZ] P.Deift, A.R.Its, X.Zhou, A Riemann-Hilbert Approach to Asymptotic Problems Arising in the Theory of Random Matrix Models and also in the Theory of Integrable Statistical Mechanics, Annals of Mathematics, v.146, n.1, 1997, pp.149-235.

[DKKZ] P.Deift, S.Kamvissis, T.Kriecherbauer, X.Zhou, The Toda Rarefaction Problem,  Comm.Pure Appl.Math v.49, n.1, 1996, pp.35-83.

[DZ1] P.Deift, X.Zhou, Long Time Asymptotics for the Autocorrelation Function of the Transverse Ising Chain at the Critical Magnetic Field,
in Singular Limits of Dispersive Waves, Plenum Press, Physics Series B320, Plenum Press, 1994.

[DZ2] P.Deift, X.Zhou, Annals of Mathematics, v.137, 1993, pp.295-368.

[DVZ] P.Deift, S.Venakides, X.Zhou, Comm.Pure Appl.Math., v.47, 1994, pp. 199-206.

[FIK] H.Frahm, A.R.Its, V.E.Korepin, An Operator-Valued Riemann-Hilbert Problem Associated with the XXX Model, preprint 95-6, IUPUI; CRM Proc.Lecture Notes, v.9, AMS 1996, pp.133-142.

[K1] S.Kamvissis, Long Time Asymptotics for the Doubly Infinite Toda Lattice, Comm.Math.Phys, v.153, n.3, 1993, pp.479-519.

[K2] S.Kamvissis, Long Time Asymptotics for  Focusing NLS with Real Spectral Singularities, Comm.Math.Phys., v.180, n.2, 1996, pp.325-342.

[KIB] V.E.Korepin, A.G.Izergin, N.M.Bogoliubov,  Quantum Inverse Scattering Method and Correlation Functions, Cambridge University Press, 1993.

[KIEU] V.E.Korepin, A.G.Izergin, F.H.L.Essler, D.B.Uglov, Correlation Function  of the Spin-1/2 XXX Antiferromagnet, Phys.Lett.A v.190, 1994, pp.182-184.

\enddocument